\documentclass[12pt]{amsart}
\usepackage{SSdefn}

\title{Structures in representation stability}
\author{Steven V Sam}
\address{Department of Mathematics, University of California, San Diego, CA}
\email{\href{mailto:ssam@ucsd.edu}{ssam@ucsd.edu}}
\urladdr{\url{http://math.ucsd.edu/~ssam/}}
\thanks{SS was supported by NSF grant DMS-1849173 and a Sloan Fellowship.}

\newcommand{\OI}{\mathbf{OI}}
\newcommand{\FI}{\mathbf{FI}}
\newcommand{\FS}{\mathbf{FS}}

\newcommand{\VIC}{\mathbf{VIC}}

\begin{document}

\maketitle

\section{Introduction}

Representation theory is applicable in many other areas of mathematics because it can be used to exploit symmetries to simplify calculations. There are many cases where the relevance is clear, such as the action of the invertible matrices on vector spaces, linear maps, tensor products, etc. via change of basis, or the action of the symmetric group by permuting coordinates or points in a space. In the examples we will discuss, we start with a sequence of objects with group actions. Next, we construct a sequence of vector spaces that are associated to the sequence of objects that carry an induced linear action of the groups. Two common examples of such groups are symmetric groups and general linear groups. At least over the field of complex numbers, the representation theory of these two groups is well-understood and allows us to group together vectors by considering irreducible decompositions. Sometimes this isn't enough though: the dimension of these vector spaces might grow fast compared to the number of isomorphism classes of irreducible representations of these groups and hence these groupings can get unwieldy.

The phrase ``representation stability'', as it will be discussed in this article, refers to two related situations and goals. The common theme is that the sequence of symmetry groups is governed by a larger algebraic structure which controls how they interact with one another. The first situation involves showing that these representations follow some predictable pattern, or stabilize in an appropriate sense. In other situations, the representations may have additional structure, such as being rings, and the explicit patterns themselves are not of interest. Instead, one may like to find bounds on invariants, such as degree of generation, or at least deduce their existence. The goal of this article is to survey a few examples of both kinds.

\section{An example of using symmetry} \label{sec:example}

Before launching into examples of representation stability, we start with an example which illustrates using symmetry to simplify calculations, and how it naturally leads to the issues studied in representation stability.

Let $V_1,\dots,V_n$ be complex finite-dimensional vector spaces and let $\bV = V_1 \otimes \cdots \otimes V_n$ denote their tensor product, a vector space (whose elements are called tensors) spanned by symbols of the form $v_1 \otimes \cdots \otimes v_n$ with $v_i \in V_i$ subject to relations to make them multilinear in the factors. The elements of the form $v_1 \otimes \cdots \otimes v_n$ are called simple tensors. We define the rank of a tensor to be the minimal $r$ such that it can be expressed as a linear combination of $r$ simple tensors. 

When $n=2$, one can identify $V_1 \otimes V_2$ with the space of linear maps from $V_1^*$ to $V_2$, and this notion of rank coincides with the usual rank of a linear map. While all mathematicians know how to compute the rank of a linear map, say using Gaussian elimination, a less traditional way is that a linear map has rank $\le r$ if and only if every $r \times r$ square submatrix (once bases are chosen) has determinant equal to $0$. This describes the locus of rank $\le r$ matrices as the zero locus of a collection of polynomials; a basic question is to find a similar description when $n>2$.\footnote{Standard caveat: the set of tensors of rank $\le r$ is generally not closed in the Zariski topology, so in what follows, we will implicitly be dealing with its Zariski closure.}

This is a largely open problem and we will focus on the case when $n=5$, $\dim V_i = 2$ for all $i$, and $r=5$. In this case, it is known that the dimension of the rank $\le 5$ locus is 2 less than the dimension of $\bV$, and so one might hope to find two polynomials whose zero locus is this set. We set off on this task in \cite{oeding-sam}: we easily found a polynomial of degree 6 which is identically 0 on this set, and experimental computation suggested that the other polynomial has degree 16.\footnote{Any multiple of the degree 6 polynomial will be identically 0 on the set of rank $\le 5$ tensors, and provides no new information. So we are implicitly discussing minimal equations, i.e., none of them are obtained from the others by multiplication and addition.} Finding such a polynomial amounts to a linear algebra problem: we first evaluate all monomials of degree 16 in $2^5=32$ variables on some set $S$ of rank $5$ tensors and organize the result into a matrix; kernel elements of this matrix give polynomials that are 0 on $S$, and if $S$ is sufficiently large and ``general'' then it gives the  desired polynomial. Practically, this is impossible: the number of such monomials is $\binom{47}{31} = 1503232609098$ (approximately 1.5 trillion).

However, we can take advantage of the fact that tensor rank is invariant under change of basis on each $V_i$, i.e., the group $G = \GL_2(\bC) \times \GL_2(\bC) \times \GL_2(\bC) \times \GL_2(\bC) \times \GL_2(\bC)$. General principles tell us that if this degree 16 polynomial exists, then we can find one which is invariant under $G$ up to scaling. This is much better since the space of such polynomials is only 1313-dimensional, but we can do even better by noting that there is also an action of the symmetric group $S_5$ which permutes the tensor factors and also leaves tensor rank unaffected. Again, such a polynomial must be invariant under this extra $S_5$-action up to scaling. All together this space is 49-dimensional, which was enough of a reduction for us to find the polynomial.

The obvious objections the reader might have: Why $n=5$? Why $\dim V_i=2$? And why $r=5$? What is so special about these parameters? Aside from the fact that the codimension is 2 in this case (codimension 1 is much easier to understand), nothing in particular is special. One would like to understand general $n$, with $\dim V_i$ general, and $r$ general. However, the information obtained in this case can be used to find equations in other cases: ``inheritance'' and ``flattening'' (which we won't make precise here) allow one to lift these equations whenever $n \ge 5$ and whenever $\dim V_i \ge 2$ (but keeping $r$ fixed). Hence one can think of these two equations as generating further equations; a better problem is to find all of the generators (say when $r=5$) under these operations rather than studying the problem one set of parameters at a time.

Even more fundamental questions are: How do we axiomatize this algebraic structure given by the operations of inheritance and flattening? As we vary over all choices of $n$ and $\dim V_i$, are the set of (minimal) equations generated by a finite number of equations under these operations? The remainder of the article is devoted to surveying other situations where a similar setup can be found. Surprisingly, the examples come from very different parts of mathematics, but the ideas needed to address the questions overlap in fundamental ways.

\section{Stability and patterns in representations} \label{sec:rep-stab}

In this section, we will consider a sequence of objects $X_0, X_1, X_2, \dots$ usually with a group action, and how we can find or study patterns that exist in numerical and linear invariants of the $X_i$. We will focus on situations where one can find some large algebraic structure that controls all of them at once and what can be deduced from this structure. We will stick to complex vector spaces for simplicity of exposition.

\subsection{A perspective on homological stability} \label{ss:hom}

Let $\bk$ be a field and let $X_n$ be the group $\GL_n(\bk)$ of $n \times n$ invertible matrices with coefficients in $\bk$. We would like to understand the $i$th homology of $X_n$, for $i$ fixed and $n$ varying. Put $M_n=\rH_i(X_n, \bC)$, i.e., the $i$th left derived functor of the functor $V \mapsto V \otimes_{X_n} \bC$. The inclusion $X_n \subset X_{n+1}$ given by $A \mapsto \begin{bmatrix} A & 0 \\ 0 & 1 \end{bmatrix}$ induces a map $M_n \to M_{n+1}$. Thus we obtain the following system of vector spaces:
\begin{displaymath}
M_0 \to M_1 \to M_2 \to \cdots
\end{displaymath}
Putting $M=\bigoplus_{n \ge 0} M_n$, we see that $M$ is a graded module over the polynomial ring in one variable $\bC[t]$, where $t$ has degree one (multiplication by $t$ on an element in $M_n$ is defined by the map $M_n \to M_{n+1}$); this exactly captures the structure we see. The most obvious question to ask is: is $M$ finitely generated as a $\bC[t]$-module? If the answer to the question is ``yes,'' then the structure theorem for finitely generated $\bC[t]$-modules tells us that $M$ decomposes as $T \oplus F$, where $T$ is a finitely generated torsion module and $F$ is a finite rank free module. This translates to a concrete statement about the original invariants: once $n$ exceeds the maximal degree in $T$, the map $\rH_i(X_n, \bC) \to \rH_i(X_{n+1}, \bC)$ is an isomorphism; that is, the homology of $\GL_n(\bk)$ stabilizes. In fact, it is known that the homology stabilizes for all fields $\bk$ and for a large class of coefficient rings besides $\bC$.

\subsection{FI-modules and cohomology of configuration spaces}

Let $Y$ be a fixed topological space and let $X_n$ be the configuration space of $n$ distinct labeled points in $Y$; thus $X_n$ is the open subset of the cartesian power $Y^n$ where the coordinates are required to be distinct. Configuration spaces appear in many places in mathematics, and it is an important problem to understand their cohomology. For example, if $Y=\bR^2$ is the plane, then $X_n$ is an Eilenberg--MacLane space for the $n$th pure braid group, and so its cohomology is the cohomology of this important group.

As in the previous example, we fix a cohomological degree $i$ and let $n$ vary. Thus put $M_n = \rH^i(X_n, \bC)$. The $n$th symmetric group acts on $X_n$ by permuting coordinates, and this induces an action on $M_n$. In particular, we are in the situation described in the introduction: we have a sequence $M_n$ of $S_n$-representations and we want to understand patterns that occur as we vary $n$. For special cases of $Y$, this can be worked out explicitly. For example, a classical calculation of Arnol'd determines these spaces when $Y=\bR^2$.

However, we want to find some intrinsic structure which might help us deduce things about a general class of manifolds. Permutations are bijective functions, and the key observation is that injective functions give extra symmetries if we consider all $X_n$ at once. For notation, set $[m] = \{1,2,\dots,m\}$ for a non-negative integer $m$. Given an injection $f \colon [m] \to [n]$, we get a map $f^* \colon X_n \to X_m$ defined by $f^*(x_1,\dots,x_n) = (x_{f(1)}, \dots, x_{f(m)})$.
Thus, after taking cohomology, we obtain a linear map $f_* \colon M_m \to M_n$. These satisfy $(f\circ g)_*=f_* \circ g_*$ for any $g \colon [n] \to [p]$. If $f$ is a permutation, then $f_* \colon M_n \to M_n$ is the action we had before, so we are extending the action of permutations on $\bigoplus_n M_n$ to an action of all injective functions ($f_*$ acts by $0$ on $M_p$ if $[p]$ is not the domain of $f$).

To formalize this, let $\FI$ (= finite injections) be the category whose objects are $[n]$ for $n \ge 0$ and whose morphisms are injective functions. Then $M$ is a functor from $\FI$ to the category of vector spaces, or more shortly, an $\FI$-module: for every $n$ we have a vector space $M_n$, and for every morphism $[m] \to [n]$ in $\FI$ we have a linear map $M_m \to M_n$, so that composition is respected. For a fixed $n$, the set of injections $[n] \to [n]$ is closed under composition and gives the symmetric group $S_n$. So, for every $\FI$-module $M$, $M_n$ is a representation of $S_n$. Hence an $\FI$-module is a sequence of $S_n$-representations together with transition maps between the different representations.

Alternatively, we say that an $\FI$-module is a representation of the category $\FI$. More generally, given a category $\cC$, a representation of $\cC$ (or $\cC$-module) is a functor from $\cC$ to the category of vector spaces.

\subsection{Finite generation}

Unfortunately, the structure of an $\FI$-module by itself is not helpful: {\it any} sequence of $S_n$-representations can be upgraded to an $\FI$-module by declaring that $f_*$ is the $0$ map whenever $f$ is not a bijection. From the example of $\bC[t]$-modules in \S\ref{ss:hom}, we see the same phenomenon: any sequence of vector spaces can be made into a graded $\bC[t]$-module by having $t$ act by $0$, so that there is no control at all over their dimensions. But, we saw that requiring finite generation avoids this problem.

So we ask this question: given the $\FI$-module $M$ coming from the space $Y$, is $M$ finitely generated, that is, are there finitely many cohomology classes that give rise to all cohomology classes by applying the $\FI$-operations and taking linear combinations? Independent of any structure theorem for $\FI$-modules, such a result is of interest since it is a kind of bound on the complexity of the cohomology classes of $X_n$ as $n$ grows. However, there is a structure theorem for finitely generated $\FI$-modules. Before stating it, we recall that irreducible complex representations of $S_n$ are indexed by integer partitions of $n$, we denote the representation corresponding to a partition $\lambda$ by $\bM_\lambda$. If $M$ is a finitely generated $\FI$-module, then the $M_n$'s are ``representation stable'' in the sense of Church--Farb; this means that there are partitions $\lambda^1, \ldots, \lambda^r$ such that for $n \gg 0$ the decomposition of $M_n$ into Specht modules (the irreducible representations of $S_n$) is $\bigoplus_{i=1}^r \bM_{\lambda^i[n]}$, where $\lambda^i[n] = (n-|\lambda^i|,\lambda^i_1, \lambda^i_2,\dots)$. A concrete consequence is that if $M$ is a finitely generated $\FI$-module, then the sequence $n \mapsto \dim M_n$ agrees with a polynomial function for $n \gg 0$.

\begin{example}
To illustrate some of these ideas in a toy example, consider the $\FI$-module $M$ where $M_n = \bC^n$ and for $f \colon [m] \to [n]$, we define $f_* \colon \bC^m \to \bC^n$ by $f_*(e_i) = e_{f(i)}$ where $e_1,e_2,\dots$ are the standard basis vectors. Then $M$ is finitely generated by one element $e_1 \in M_1$. When $n \ge 2$, $M_n$ decomposes into two irreducible representations: the subspace $\{(x_1,\dots,x_n) \mid x_1+ \cdots+x_n =0\}$ and the line spanned by $e_1+\cdots+e_n$. The first is $\bM_{(n-1,1)}$ while the latter is $\bM_{(n)}$. In this case, the partitions are $\lambda^1=(1)$ and $\lambda^2=\emptyset$.
\end{example}

One of the first main results about $\FI$-modules, due to Church, Ellenberg, and Farb \cite{CEF}, is that the $\FI$-module $M$ associated to $Y$ is finitely generated when it is a manifold with some mild restrictions. 
There is a wealth of literature surrounding $\FI$-modules and their applications, see \cite{farb} for some further references.

So it is desirable to establish finite generation of certain representations of categories. In most applications, this is done in two steps:
\begin{enumerate}
\item Show that representations of the category in question have the noetherian property: any subrepresentation of a finitely generated representation is again finitely generated. For $\FI$-modules, this was the main result of \cite{CEF}.
\item Use the noetherian property to prove finite generation of the specific representations in question. (For example, Church, Ellenberg, and Farb use that $\rH^i(X_n)$ can be computed by a spectral sequence of $\FI$-modules, and each module on the initial page is easily shown to be finitely generated.)
\end{enumerate}
In joint work with Snowden \cite{catgb}, we construct a general theory that establishes the noetherian property for many of the categories of interest in representation stability.

\subsection{Other examples}

We now give a brief survey of some of the other categories which have come up in representation stability. In the interest of saving space, we will not motivate the definitions or explain consequences of finite generation.

There are various ways to enlarge the category $\FI$ in order to deduce stronger properties about examples such as cohomology of configuration spaces. One such example is the category of non-commutative finite sets: a morphism $n \to m$ is an ordinary set function $f \colon [n] \to [m]$ together with a choice of total ordering on each fiber of $f$. This is used in \cite{EWG} to study configuration spaces of manifolds with a nowhere vanishing vector field.

One can define $q$-analogues of $\FI$-modules: given a finite field $\bF_q$, we define a category $\VIC(\bF_q)$ whose objects are non-negative integers, and a morphism $m \to n$ is an injective $\bF_q$-linear map between vector spaces $f\colon \bF_q^m \to \bF_q^n$ together with a choice of complementary subspace $C \subseteq \bF_q^n$ to the image $f(\bF_q^m)$. More generally, $\bF_q$ can be replaced by other finite commutative rings like $\bZ/\ell$. There is also a symplectic variant ${\bf SI}(\bF_q)$ where the linear maps are required to be compatible with the standard symplectic form on $\bF_q^{2m}$ and $\bF_q^{2n}$ (the complement is then chosen to be the orthogonal complement with respect to this form). These examples were studied in \cite{putman-sam} in connection with the homology of congruence subgroups, which are kernels of the following kinds of homomorphisms: $\GL_n(\bZ) \to \GL_n(\bZ/\ell)$, $\Sp_{2n}(\bZ) \to \Sp_{2n}(\bZ/\ell)$, maps from the mapping class group of a surface of genus $g$ to $\Sp_{2g}(\bZ/\ell)$, and maps from the automorphism group of a free group on $n$ generators to $\GL_n(\bZ/\ell)$. Each of these examples has the structure of a finitely generated representation of one of these categories.

The moduli space $\cM_{g,n}$ of genus $g$ Riemann surfaces with $n$ marked points has an important compactification $\ol{\cM}_{g,n}$ which was constructed by Deligne and Mumford. Its homology and cohomology both admit the structure of an $\FI$-module but fail to be finitely generated. However, a variant can be used: define $\FS^\op$ to be the category whose objects are non-negative integers and such that a morphism $m \to n$ is a surjection $[n] \to [m]$ (the op refers to opposite direction). Tosteson showed in \cite{tosteson} that the homology of $\ol{\cM}_{g,n}$ carries a representation of $\FS^\op$ which is finitely generated.

Let $\OI$ be the subcategory of $\FI$ where injections are required to be order-preserving. Every $\FI$-module is automatically an $\OI$-module by restricting the action. However, there are also naturally occurring examples of $\OI$-modules that don't come from $\FI$. One such involves the homology of groups of upper-triangular matrices (the ordering on the basis elements becomes important) which are shown to be finitely generated in \cite{PSS}.

\section{Existence of uniform bounds}

The functorial perspective from \S\ref{sec:rep-stab} was initially driven by topological examples in the works of Church, Ellenberg, Farb, Putman, etc. In parallel, other forms of representation stability were being used in commutative algebra and algebraic geometry, as we hinted at in \S\ref{sec:example}. Many of these results are of a much different character, though the underlying theme of proving a noetherianity result remains the same. For a survey on this side of the story, we also refer to \cite{draisma-notes}.

\subsection{Equivariant rings and modules}

Let $R$ be a ring on which a group $G$ acts by ring automorphisms. An \emph{equivariant $R$-module} is an $R$-module $M$ equipped with an action of $G$ that is compatible with its action on $R$, in the sense that $g(ax)=(ga)(gx)$ for all $g \in G$, $a \in R$, $x \in M$. Equivariant modules are ubiquitous, and the novel aspect in representation stability is that the objects involved tend to be ``large''. Some examples:
\begin{itemize}
\item Take $R$ to be the infinite variable polynomial ring $\bC[x_1,x_2,\ldots]$ and $G$ to be the infinite symmetric group $S_{\infty}$, acting by permuting the variables.
\item Take $R=\bC[x_1,x_2,\ldots]$ as above and $G$ to be the infinite general linear group $\GL_{\infty}$, acting by linear substitutions in the variables.
\item Take $R=\bC[x_{i,j}]_{i,j \ge 1}$ with $x_{i,j}=x_{j,i}$, and take $G=\GL_{\infty}$. The variables are the entries of an infinite symmetric matrix $A=(x_{i,j})$, and $g x_{i,j}$ is the $(i,j)$ entry of $gAg^T$.
\item Generalizing the previous two: let $V$ be a representation of $G=\GL_\infty$ and let $G$ act on the symmetric algebra $\Sym(V)$ (i.e., picking a basis for $V$ identifies $\Sym(V)$ with the polynomial ring with those basis elements as variables). In the first case, $V = \bC^\infty$, thought of as column vectors, and in the second, $V = \Sym^2(\bC^\infty)$ is the space of quadratic polynomials in $x_1,x_2,\dots$.
\end{itemize}

Why infinitely many variables? In applications, we would really be interested in finitely many variables, like $\bC[x_1,\dots,x_n]$ under the action of a smaller group such as the $n$th symmetric group $S_n$. Much like in the previous section, our calculation or object of interest varies with $n$ and so we get a sequence. This can sometimes be rephrased as a single object in the case when $n \to \infty$. Hence this offers a different perspective: prove properties about a large algebraic structure as opposed to a structure which governs sequences of representations.

In each case, for applications, we are interested in whether an analogue of a noetherian property holds. For example, we can ask whether ideals closed under the group action are finitely generated up to this action, i.e., given such an ideal $I$, we can ask if there are $f_1,\dots,f_r \in I$ such that every $f \in I$ is a linear combination of elements of the form $g \cdot f_i$ where $g \in G$. This property is known to hold for the first 3 examples, but not in the level of generality of the fourth one. For an application of the first example in algebraic statistics, we point to \cite{hillar-sullivant}.

There is a (seemingly) more general question to ask: if $M$ is a finitely generated equivariant $R$-module, are all equivariant submodules of $M$ also finitely generated? In standard commutative algebra, this property holds once we know that it holds for ideals, but we have thus far found no such formal implication which covers  these cases. Again, this stronger property holds for the first 3 examples, but is unknown for the fourth example.\footnote{In the second example, the category of modules with a polynomial action of $\GL_\infty$ is equivalent to the category of $\FI$-modules, see \cite[Proposition 1.3.5]{symc1}.}

\subsection{Topological noetherianity}

On the other hand, there is a weaker property we can ask for. In algebraic geometry, an ideal $I$ in a polynomial ring $\bC[x_1,\dots,x_n]$ (we allow $n=\infty$) corresponds to an algebraic set: this is the set of points in $\bC^n$ which give 0 when substituted into any polynomial in $I$. Finite generation of the ideal implies that there is a finite list of polynomials $f_1,\dots,f_r$ so that membership in the corresponding algebraic set can be tested by evaluating just these $r$ polynomials at the point. When $n$ is finite, every ideal is finitely generated by the Hilbert basis theorem. This can be rephrased as saying that $\bC[x_1,\dots,x_n]$ is a noetherian ring. This implies that $\bC^n$ is {\it topologically noetherian}, i.e., testing membership in any algebraic set can be done with a finite list of polynomials. Note that being topologically noetherian says less than the ring being noetherian since different ideals can give the same algebraic set.

In our equivariant infinite-dimensional context, $\bC^\infty$ carries an action of the group $G$ that acts on $\bC[x_1,x_2,\dots]$ and so we will be interested in algebraic sets closed under the $G$-action. Now we can ask if the space is topologically $G$-noetherian: for every equivariant algebraic set, we want a finite list of polynomials $f_1,\dots,f_r$ so that membership of $x$ can be tested by deciding if $(g \cdot f_i)(x)=0$ for all $g \in G$ and $i=1,\dots,r$. The fourth example above is known by work of Draisma \cite{draisma} to be topologically noetherian for a large class of $V$: the polynomial representations. Two particularly interesting classes of polynomial representations are:
\begin{itemize}
\item Points of the symmetric power $\Sym^d \bC^\infty$ are degree $d$ homogeneous polynomials in $x_1,x_2,\dots$, and hence when $V = \bigoplus_{i=1}^r \Sym^{d_i} \bC^\infty$, points parametrize tuples of homogeneous polynomials $(f_1,\dots,f_r)$ with $\deg(f_i)=d_i$, which can be used to parametrize ideals in any polynomial ring in finitely many variables. This perspective gives one a way to study invariants of ideals and, in particular, stabilization of such invariants in families when the degrees of the generators of our ideal are fixed in advance. 
Projective dimension is an important example of such an invariant, which connects this circle of ideas with the work on Stillman's conjecture \cite{AH, bigpoly}.
  
\item Points of the exterior power $\bigwedge^d \bC^\infty$ represent degree $d$ skew-commutative polynomials in $x_1,x_2,\dots$. The BGG correspondence allows one to study cohomology of sheaves on projective space in terms of linear algebra computations with skew-commutative polynomials. 
This perspective allows one to prove stabilization properties of cohomology of sheaves in various kinds of families, see \cite{cohomology}.
\end{itemize}

\subsection{Bounded rank tensors}

In this last part, we return to the questions raised in \S\ref{sec:example}.

First, the notions of inheritance and flattening are axiomatized by $\Delta$-modules in the sense of \cite{delta-mod}. These are functors from a category $\Delta$ (whose definition is too involved to give here) to the category of vector spaces. The objects of $\Delta$ are finite tuples of finite-dimensional vector spaces, and the assignment of a tuple to the space of degree $d$ polynomials that are $0$ on the locus of rank $\le r$ tensors is an example of a $\Delta$-module. In fact, it is finitely generated, which answers the question posed earlier, with the caveat that both $r$ and $d$ must be fixed.

If we want also to allow $d$ to vary, then we need a more sophisticated algebraic structure. It is currently unknown whether it is possible to find such a structure which acts on all vanishing polynomials in a finitely generated way. However, Draisma and Kuttler \cite{draismakuttler} give a topological version by proving the following statement: if we fix $r$, there is a constant $C(r)$, such that for every tuple of vector spaces, there exists a finite collection of polynomials of degree $\le C(r)$ which can be used to test membership in the rank $\le r$ locus.\footnote{The subtle difference here is about finding set-theoretic equations versus finding all equations, or more technically, about finding all generators for an ideal versus finding enough polynomials that generate some ideal with the same radical.} To prove this statement, Draisma and Kuttler take an appropriate limit of the tensor product $V_1 \otimes \cdots \otimes V_n$ (as $n \to \infty$ and allowing $\dim V_i$ to vary) to get an infinite-dimensional topological space with an action of a group $G$, and prove that the space is topologically $G$-noetherian.

Finally, there are several important variations of this problem with different answers given. First, instead of tensor products $V_1 \otimes \cdots \otimes V_n$, we can consider exterior powers $\bigwedge^n V$ and symmetric powers $\Sym^n V$. Simple tensors are, respectively, those of the form $v_1 \wedge \cdots \wedge v_n$ and $v^n$, so that we have analogues of tensor rank. The analogue of Draisma and Kuttler's result for $\bigwedge^n V$ is proven by Draisma and Eggermont \cite{draismaeggermont}. An answer about the full ideal (not just a topological statement) in the spirit of the original question is proven by Laudone \cite{laudone} for $\bigwedge^n V$ and by myself \cite{sam} for $\Sym^n V$.


\begin{thebibliography}{CEFN}

\bibitem[AH]{AH} Tigran Ananyan, Melvin Hochster, Small subalgebras of polynomial rings and Stillman's conjecture, \arxiv{1610.09268v3}.
  
\bibitem[CEF]{CEF} Thomas Church, Jordan Ellenberg, Benson Farb, FI-modules and stability for representations of symmetric groups, {\it Duke Math. J.} {\bf 164} (2015), no.~9, 1833--1910, \arxiv{1204.4533v3}.

\bibitem[Dr1]{draisma-notes} Jan Draisma, Noetherianity up to symmetry, {\it Combinatorial algebraic geometry}, Lecture Notes in Math. {\bf 2108}, Springer, 2014, \arxiv{1310.1705v2}.

\bibitem[Dr2]{draisma} Jan Draisma, Topological noetherianity of polynomial functors, {\it J. Amer. Math. Soc.} {\bf 32}, no.~3 (2019), 691--707, \arxiv{1705.01419v4}.

\bibitem[DE]{draismaeggermont} Jan Draisma, Rob H. Eggermont, Pl\"ucker varieties and higher secants of Sato's Grassmannian, {\it J. Reine Angew. Math.}, to appear, \arxiv{1402.1667v3}.

\bibitem[DK]{draismakuttler} Jan Draisma, Jochen Kuttler, Bounded-rank tensors are defined in bounded degree, {\it Duke Math. J.} {\bf 163} (2014), no.~1, 35--63, \arxiv{1103.5336v2}.


\bibitem[EWG]{EWG} Jordan S. Ellenberg, John D. Wiltshire-Gordon, Algebraic structures on cohomology of configuration spaces of manifolds with flows, \arxiv{1508.02430v2}.
  
\bibitem[ESS1]{bigpoly} Daniel Erman, Steven V Sam, Andrew Snowden, Big polynomial rings and Stillman's conjecture, {\it Invent. Math.}, to appear, \arxiv{1801.09852v4}.


\bibitem[ESS2]{cohomology} Daniel Erman, Steven V Sam, Andrew Snowden, Stillman uniformity for cohomology of sheaves, \arxiv{1906.10870v1}.

\bibitem[Fa]{farb} Benson Farb, Representation stability, {\it Proceedings of the International Congress of Mathematicians. Volume II}, 1173--1196, \arxiv{1404.4065v1}.
  
\bibitem[HS]{hillar-sullivant} Christopher J. Hillar, Seth Sullivant, Finite Gr\"obner bases in infinite dimensional polynomial rings and applications, {\it Adv. Math.} {\bf 229} (2012), no.~1, 1--25.

\bibitem[La]{laudone} Robert P. Laudone, Syzygies of secant ideals of Pl\"ucker-embedded Grassmannians are generated in bounded degree, \arxiv{1803.04259v2}.

\bibitem[OS]{oeding-sam} Luke Oeding, Steven V Sam, Equations for the fifth secant variety of Segre products of projective spaces, {\it Exp. Math.} {\bf 25} (2016), no.~1, 94--99, \arxiv{1502.00203v2}.

        
\bibitem[PS]{putman-sam} Andrew Putman, Steven~V Sam, Representation stability and finite linear groups, {\it Duke Math. J.} {\bf 166} (2017), no.~13, 2521--2598, \arxiv{1408.3694v3}.

\bibitem[PSS]{PSS} Andrew Putman, Steven V Sam, Andrew Snowden, Stability in the homology of unipotent groups, \arxiv{1711.11080v3}.

\bibitem[SS1]{symc1} Steven V Sam, Andrew Snowden, GL-equivariant modules over polynomial rings in infinitely many variables, {\it Trans. Amer. Math. Soc.} {\bf 368} (2016), 1097--1158, \arxiv{1206.2233v3}.
  
\bibitem[SS2]{catgb} Steven~V Sam, Andrew Snowden, Gr\"obner methods for representations of combinatorial categories, {\it J. Amer. Math. Soc.} {\bf 30} (2017), no.~1, 159--203, \arxiv{1409.1670v3}.

\bibitem[Sa]{sam} Steven~V Sam, Ideals of bounded rank symmetric tensors are generated in bounded degree, {\it Invent. Math.} {\bf 207} (2017), no.~1, 1--21, \arxiv{1510.04904v2}.

\bibitem[Sn]{delta-mod} Andrew Snowden, Syzygies of Segre embeddings and $\Delta$-modules, {\it Duke Math. J.} {\bf 162} (2013), no.~2, 225--277, \arxiv{1006.5248v4}.

  
\bibitem[To]{tosteson} Philip Tosteson, Stability in the homology of Deligne--Mumford compactifications, \arxiv{1801.03894v3}.


\end{thebibliography}
\end{document}